\documentclass[12pt]{amsart}

\usepackage{amssymb,amsmath,amsthm}
\usepackage[T1]{fontenc} % T1 for the listings' straight quotes (upquote); OT1 has no \textquotedbl
\usepackage{lmodern}
\usepackage{fullpage}
\usepackage{url}
\usepackage{color}
\usepackage{multirow}
\usepackage[table]{xcolor}
\usepackage{textcomp,listings}
\usepackage{placeins} %for \FloatBarrier
\usepackage{graphicx} %scalebox

\definecolor{green}{rgb}{0,0.5,0}
\definecolor{dkgreen}{rgb}{0,0.6,0}
\definecolor{gray}{rgb}{0.5,0.5,0.5}
\definecolor{mauve}{rgb}{0.58,0,0.82}
\definecolor{questcolor}{rgb}{0.9, 0.9, 0.9}
\lstnewenvironment{python}[1][]{
\lstset{
  frame=single,                   % adds a frame around the code
  language=python,                          % the language of the code
  basicstyle=\ttfamily\small, % the size of the fonts that are used for the code
  numbers=left,                             % where to put the line-numbers
  numberstyle=\scriptsize\color{black},  % the style that is used for the line-numbers
  stepnumber=1,                   % the step between two line-numbers. If it's 1, each line
  numbersep=8pt,                  % how far the line-numbers are from the code
  backgroundcolor=\color{white},      % choose the background color.
  stringstyle=\color{red},
  showspaces=false,               % show spaces adding particular underscores
  showstringspaces=false,         % underline spaces within strings
  showtabs=false,                 % show tabs within strings adding particular underscores
  tabsize=2,                      % sets default tabsize to 2 spaces
  captionpos=b,                   % sets the caption-position to bottom
  breaklines=true,                % sets automatic line breaking
  breakatwhitespace=false,        % sets if automatic breaks should only happen at whitespace
  rulecolor=\color{black},        % if not set, the frame-color may be changed on line-breaks within not-black text (e.g. commens (green here))
  framexleftmargin=1mm, framextopmargin=1mm, frame=shadowbox, %rulesepcolor=\color{blue},#1
  alsoletter={1234567890},
  otherkeywords={\ , \}, \{},
  keywordstyle=\color{blue},       % keyword style
  stringstyle=\color{mauve},         % string literal style
  emph={access,and,break,class,continue,def,del,elif ,else,%
  except,exec,finally,for,from,global,if,import,in,i s,%
  lambda,not,or,pass,print,raise,return,try,while},
  emphstyle=\color{black}\bfseries,
  emph={[2]True, False, None, self},
  emphstyle=[2]\color{blue},
  emph={[3]from, import, as},
  emphstyle=[3]\color{blue},
  upquote=true,
  morecomment=[s]{"""}{"""},
  commentstyle=\color{dkgreen},       % comment style
  emph={[4]1, 2, 3, 4, 5, 6, 7, 8, 9, 0},
  emphstyle=[4]\color{blue},
  literate={:}{{\textcolor{blue}:}}{1}%
  {=}{{\textcolor{blue}=}}{1}%
  {-}{{\textcolor{blue}-}}{1}%
  {+}{{\textcolor{blue}+}}{1}%
  {!}{{\textcolor{blue}!}}{1}%
  {(}{{\textcolor{blue}(}}{1}%
  {)}{{\textcolor{blue})}}{1}%
  {[}{{\textcolor{blue}[}}{1}%
  {]}{{\textcolor{blue}]}}{1}%
  {<}{{\textcolor{blue}<}}{1}%
  {>}{{\textcolor{blue}>}}{1},%
}}{}

\newcommand{\hhat}{{\hat h}}

\providecommand{\abs}[1]{\left\lvert#1\right\rvert}

\def\bsigma{\boldsymbol{\sigma}}

\def\A{\mathbb{A}}
\def\C{\mathbb{C}}

\def\P{\mathbb{P}}
\def\Q{\mathbb{Q}}

\def\Z{\mathbb{Z}}

\def\calM{\mathcal{M}}

\def\calO{\mathcal{O}}

\def\Kbar{\overline{K}}
\def\Qbar{\overline{\Q}}

\DeclareMathOperator{\Aut}{Aut}

\DeclareMathOperator{\Hom}{Hom}
\DeclareMathOperator{\Per}{Per}
\DeclareMathOperator{\PGL}{PGL}
\DeclareMathOperator{\Pre}{Pre}

\newcommand{\col}{\,{:}\,}

\theoremstyle{plain}% default
\newtheorem{thm}{Theorem}

\theoremstyle{definition}
\newtheorem{defn}[thm]{Definition}
\newtheorem{conj}{Conjecture}
\newtheorem{exmp}[thm]{Example}
\newtheorem{rem}{Remark}
\theoremstyle{remark}

\title[Dynabase]{Dynabase: A Database of Extreme Examples in Arithmetic Dynamics}

\author{Benjamin Hutz}
\address{Department of Mathematics and Statistics,
         Saint Louis University,
         St.~Louis, MO 63103}
\email{benjamin.hutz@slu.edu}

\begin{document}



% TODO: funding and student acknowledgments; check where the target journal wants the AI-use disclosure.
\thanks{The author acknowledges the use of Claude Opus 5 and 5.5 in the creation of the GitHub Pages site hosting Dynabase as well as in copyediting the manuscript.}

%\date{\today}

    \begin{abstract}
        We introduce Dynabase, a database cataloging extreme examples in arithmetic dynamics related to several major open conjectures: the Morton--Silverman uniform boundedness conjecture and the dynamical analogue of Lang's conjecture on canonical heights, as well as data on postcritically finite maps and maps with nontrivial automorphisms. This database collects the known examples of extreme behavior for these problems and tracks the current status of subproblems based on degree and field of definition.
    \end{abstract}

\maketitle

Algebraic dynamics is a field rich in conjectured behaviors and expected uniform constants. However, only a handful of specific cases have been resolved or computed explicitly, and these results are scattered throughout the research literature. The goal of Dynabase is to drive further work on these problems through two mechanisms:
\begin{enumerate}
    \item Collect and organize the known examples into a single browsable and downloadable location.
    \item Track the status of various subproblems by degree and field of definition.
\end{enumerate}
The examples in Dynabase are related to four main problems:
\begin{enumerate}
    \item Rational preperiodic points
    \item Postcritically finite maps
    \item Automorphism groups
    \item Small canonical height ratios
\end{enumerate}
Each problem is introduced in detail in Section \ref{sect_background}.

Dynabase is a public GitHub project which welcomes contributors. It can be thought of as having three main parts:
\begin{enumerate}
    \item Collecting the examples
    \item Computing and storing their properties
    \item A web front-end for browsing and exploring
\end{enumerate}
\textit{Collecting the examples:} The collection of examples comes from publications and personal communications. There is also a mechanism for future submissions on the website as well as the ability to submit pull requests through GitHub. The goal is not to collect every map with extreme properties; the goal is to collect all the possible extreme properties. For example, a quadratic polynomial map defined over $\Q$ is conjectured to have at most $9$ $\Q$-rational preperiodic points, counting the point at infinity in $\P^1(\Q)$ \cite{Poonen}. Poonen proved this bound provided there are no rational periodic points of minimal period greater than $3$; periods $4$ and $5$ are ruled out \cite{FPS, Morton5}, and period $6$ is ruled out assuming the Birch and Swinnerton-Dyer conjecture for a certain Jacobian \cite{Stoll3}. However, there are infinitely many distinct such maps that exhibit this property. Only one such map per graph structure is stored in the database because the other maps exhibit the same extreme behavior. There is currently no clear criterion as to which map with a given property should be included; any new map is merely checked for being a duplicate before being added.

\textit{Computing and storing their properties:} Computing the properties is done mainly through SageMath. A significant amount of work is done to normalize the map as well as the field of definition. Two maps in the same conjugacy class are typically treated as the same; however, note that maps which are rational twists may have different properties with respect to rationality questions (such as the number of rational preperiodic points). In that situation, it is possible for two rational twists to both occur in the database. The data is stored in a local PostgreSQL database and can be populated simply by running a SageMath script (after the needed applications are installed). The majority of the data can be populated in about 8 hours of computation on a standard desktop computer, but a few of the properties/maps required some additional computation time.

\textit{A web front-end for browsing and exploring:}
Currently, the web front-end is a set of static pages which include
\begin{enumerate}
    \item General background on dynamics and the problems
    \item Summary of data currently included
    \item Specific data pages for each subproblem
    \item Summary of extreme examples
    \item Summary of current status of the subproblems.
\end{enumerate}
There are several future directions this project can go.
\begin{enumerate}
    \item Continue to add data as more subproblems get addressed.
        \begin{enumerate}
            \item Additional degrees and fields
            \item Higher dimensional maps
        \end{enumerate}
    \item Add additional subproblems.
    \item Add problems over finite fields.
\end{enumerate}

Dynabase follows the model of databases such as the LMFDB \cite{LMFDB} and the databases of number fields \cite{Jones4} and of genus 2 curves \cite{Booker}, which utilize systematically computed and normalized data, with stable labels as the basis for driving future research. The emphasis on extreme examples is similar to Dujella's tables of elliptic curves of high rank with prescribed torsion \cite{Dujella}, which list the current records rather than all known curves. In arithmetic dynamics, a number of publications have explored these problems from a computational perspective. For example, the problem of rational preperiodic points of quadratic polynomials over $\Q$ and quadratic fields is examined in \cite{JDoyle, DFK, Poonen}. The problem of small canonical heights is examined in \cite{Benedetto8, Benedetto4, Hutz29}. The problem of postcritically finite maps is examined in \cite{Anderson,Ingram2,Lukas}. The problem of identifying maps and computing maps with automorphisms is examined in \cite{Hutz16, Hutz25, Miasnikov}. These are just a few examples of articles addressing these problems and we already see that the results are spread across many papers, each with its own normalizations and conventions. Dynabase collects this data in one place, in a common normalization, and tracks how it relates to the open problems.

The organization of this article is as follows. Section \ref{sect_background} gives background on arithmetic dynamics and describes each of the four problems in detail. Section \ref{sect_scope} describes what Dynabase stores and which examples count as distinct, and Section \ref{sect_computation} describes how the data is computed. Section \ref{sect_labels} proposes a labelling system for dynamical systems. Section \ref{sect_snapshot} summarizes the current data release. Section \ref{sect_corrections} lists errors in the literature found while compiling the data. Section \ref{sect_using} describes the website and the export to SageMath, and Section \ref{sect_architecture} describes the software and how to contribute. Section \ref{sect_AI} describes the AI use in this project and Section \ref{sect_future} discusses future directions.

\section{Dynamics Background} \label{sect_background}
We restrict to the case of number fields, as those are the problems currently addressed in Dynabase, but the definitions are valid in much more general settings. Let $K$ be a number field and $\P^1$ the projective line. An endomorphism $f: \P^1 \to \P^1$ of degree $d \geq 2$ is defined by a pair of homogeneous polynomials of degree $d$ with no common factors. We denote by $\Hom_d$ the set of all endomorphisms of degree $d$ on $\P^1$ and by $\Hom_d(K)$ those defined over $K$. We can iterate $f$ to create a dynamical system. Specifically, we define the $n$-th iterate of $f$ recursively as $f^0 = \mathrm{id}$ and $f^n = f \circ f^{n-1}$ for $n \geq 1$. Since $f$ is defined by polynomial equations, we consider this the study of algebraic dynamical systems. When the questions asked have a number theoretic flavor to them (such as rationality and field of definition questions), we term this arithmetic dynamics. All current problems in Dynabase can be classified as arithmetic dynamical problems.

One class of problems surrounds the notion of preperiodic points, which are points with finite forward orbit. Specifically, we define
\begin{equation*}
    \calO_f(P) = (P, f(P), f^2(P), \dots)
\end{equation*}
as the \emph{forward orbit} of $P$ by $f$. If the set $\{f^n(P) \col n \geq 0\}$ is finite, we say that $P$ is \emph{preperiodic for $f$}. We can make the further distinction that $P$ is \emph{periodic with period $n$} if $f^n(P) = P$ for some $n \geq 1$, and $P$ is \emph{preperiodic with period $(m,n)$} if $f^m(P)$ is periodic with period $n$. The value $m$ is called the \emph{tail} of $P$. The minimal such values of $n$ and $m$ are called the \emph{minimal period} and the \emph{minimal tail} of $P$. We denote by $\Per_n(f)$ the set of points of (not necessarily minimal) period $n$ and by $\Pre(f,K)$ the set of $K$-rational preperiodic points of $f$. Some questions that arise are
\begin{enumerate}
    \item What is the maximum of $\#\Pre(f,K)$ over all $f \in \Hom_d(K)$?
    \item What is the longest possible minimal period of a $K$-rational periodic point for some $f \in \Hom_d(K)$?
    \item What is the longest possible minimal tail length of a $K$-rational preperiodic point for some $f \in \Hom_d(K)$?
\end{enumerate}
These questions are related to the Morton--Silverman uniform boundedness conjecture discussed in Section \ref{sect_preper}.

We next look at the critical points and postcritically finite maps.
\begin{defn}
    Let $f\in \Hom_d$ for $d \geq 2$. The \emph{critical points of $f$} (or \emph{ramification locus}) are the points $P$ where $f$ is not injective in any neighborhood of $P$. Specifically, writing $f = (F:G)$ for homogeneous polynomials $F$ and $G$, these are the points where the Jacobian determinant $F_xG_y - F_yG_x$ vanishes. Counted with multiplicity, there are $2d-2$ critical points.

    We say $f$ is \emph{postcritically finite} if all of its critical points are preperiodic.
\end{defn}
Postcritically finite maps are discussed in Section \ref{sect_pcf}.

The remaining two problems involve the conjugation action. The projective linear group $\PGL_2$ is the group of automorphisms of $\P^1$ and has a natural action on $\Hom_d$ that preserves dynamical properties (but not necessarily arithmetic properties). Let $f \in \Hom_d$ and $\alpha \in \PGL_2$; then the \emph{conjugation} of $f$ by $\alpha$ is defined as
\begin{equation*}
    f^{\alpha} = \alpha^{-1} \circ f \circ \alpha.
\end{equation*}
This preserves the dynamical properties since the iterate of the conjugate is the conjugate of the iterate
\begin{equation*}
    (f^{\alpha})^n = (f^n)^{\alpha}.
\end{equation*}
Taking the quotient of $\Hom_d$ by this conjugation action gives the moduli space of dynamical systems of degree $d$, $\calM_d$. That $\calM_d$ is well defined as a variety is proven by Silverman for $\P^1$ \cite{Silverman9} and Petsche, Szpiro, and Tepper \cite{petsche} or Levy \cite{Levy} for higher dimensions. Elements of $\calM_d$ are termed \emph{conjugacy classes}, and we say $f$ and $g$ are conjugate if there exists an element $\alpha \in \PGL_2$ so that $f^{\alpha} = g$. If $f,g$ are defined over $K$ and are conjugate (over $\Kbar$), but not by any $\alpha \in \PGL_2(K)$, then we say that $f$ and $g$ are \emph{rational twists over $K$}.

The \emph{automorphism group of $f$} is defined as
\begin{equation*}
    \Aut(f) = \{ \alpha \col f^{\alpha} = f\} \subset \PGL_2.
\end{equation*}
$\Aut(f)$ is a group under composition and always contains at least the identity element. Further, $\Aut(f)$ is a finite group since, for example, it must permute the set $\Per_n(f)$ for every $n \geq 1$. The finite subgroups of $\PGL_2(\C)$ were classified classically by Klein and presented in modern notation by Silverman \cite{Silverman12}. Throughout, $\Aut(f)$ denotes the automorphism group over $\Kbar$, even when $f$ is defined over $K$. The two main questions surrounding automorphisms are the following:
\begin{enumerate}
    \item Which maps have a nontrivial automorphism group?
    \item Which finite subgroups of $\PGL_2(\Kbar)$ are realizable as $\Aut(f)$ for some $f$ defined over $K$?
\end{enumerate}
These problems are discussed in Section \ref{sect_aut}.

The last problem addressed in Dynabase involves the canonical height. Let $M_K$ be the set of places of the number field $K$. We define
\begin{defn}
    The \emph{(relative) multiplicative height} of $a \in K$ is defined as
    \begin{equation*}
        H_K(a) = \prod_{v \in M_K} \max(\abs{a}_v,1)^{n_v},
    \end{equation*}
    where $n_v = [K_v:\Q_v]$ is the local degree. More generally, we define
    the \emph{(relative) multiplicative height} of $P=(P_0:P_1) \in \P^1(K)$ as
    \begin{equation*}
        H_K(P) = \prod_{v \in M_K}\max(\abs{P_0}_v, \abs{P_1}_v)^{n_v}.
    \end{equation*}
    The \emph{(relative) logarithmic height} of $P \in \P^1(K)$ is defined as
    \begin{equation*}
        h_K(P) = \log(H_K(P)),
    \end{equation*}
    and the \emph{(absolute) logarithmic height} of $P$ is $h(P) = \frac{1}{[K:\Q]}h_K(P)$, which does not depend on the choice of the field $K$ containing the coordinates of $P$.
\end{defn}
It is known that the height changes in a controlled manner under iteration \cite[Theorem 3.11]{Silverman10}:
\begin{equation*}
    \abs{h(f(P)) - dh(P)} \leq C \quad \text{for all } P \in \P^1(\Kbar),
\end{equation*}
where the constant $C$ depends only on $f$ and not on $P$. We can also define a canonical height
\begin{equation*}
    \hhat_f(P) = \lim_{n \to \infty} \frac{h(f^n(P))}{d^n}.
\end{equation*}
This limit exists and provides a function (the \emph{canonical height with respect to $f$}) which satisfies
\begin{equation*}
    \hhat_f(f(P)) = d\hhat_f(P).
\end{equation*}
This forces any preperiodic point to have canonical height $0$ and, for number fields, the converse is also true: any point with canonical height $0$ is a preperiodic point. When combined with the existence of a constant $C$ (depending only on $f$) so that
\begin{equation*}
    \abs{\hhat_f(P) - h(P)} \leq C,
\end{equation*}
and Northcott's theorem that the set of preperiodic points is a set of bounded height, we see that there must be some point in $\P^1(K)$ with smallest positive canonical height. A dynamical version of Lang's conjecture (see Section \ref{sect_height}) predicts how small this smallest canonical height can be.

We now describe each of these problems in more detail.

\subsection{Rational Preperiodic Points} \label{sect_preper}

    In 1950 Northcott proved that for any given map the set of preperiodic points is a set of bounded height and hence has only finitely many points defined over any given number field \cite{Northcott}. In 1994, Morton and Silverman went one step further and conjectured that the number of rational preperiodic points should be bounded independently of $f$. In particular, the upper bound should only depend on the degree $d$ of the map and the degree $D$ of the number field. They stated the conjecture for $\P^N$ (where the constant also depends on $N$), but we restate it for $\P^1$ only.
    \begin{conj}[\cite{Silverman7}] \label{conj1}
        For any integers $d \geq 2$ and $D \geq 1$, there exists a constant $C(d,D)$ so that for every number field $K$ with $[K:\Q] \leq D$ and every endomorphism $f \in \Hom_d(K)$, we have
        \begin{equation*}
            \#\Pre(f,K) \leq C(d,D).
        \end{equation*}
    \end{conj}
    For polynomial maps, Looper proved this conjecture assuming the $abcd$-conjecture \cite{Looper3}. Conjecture \ref{conj1} would have strong consequences beyond dynamics: Fakhruddin showed that the conjecture for $\P^N$ implies the uniform boundedness of torsion points on abelian varieties \cite{Fakhruddin}. Unconditionally, it is not known in any case, even for quadratic polynomials over $\Q$. For quadratic polynomials, Poonen classified the possible graphs over $\Q$ assuming there are no rational cycles of length greater than $3$ \cite{Poonen}, and Doyle, Faber and Krumm \cite{DFK} and Doyle \cite{JDoyle,JDoyle3} extended the classification to quadratic fields, conjecturally. In higher degrees, the known examples come from searches \cite{Benedetto4,Doyle-Hyde,Hutz29}. If this conjecture were true, it would imply the existence of a number of other uniform constants.
    \begin{itemize}
        \item There would be a uniform upper bound on the number of distinct directed graphs of rational preperiodic points for $f \in \Hom_d(K)$ defined over a field $K$ with $[K:\Q] \leq D$.
        \item There would be a uniform upper bound on the minimal period of a $K$-rational periodic point.
        \item There would be a uniform upper bound on the minimal tail length of a $K$-rational preperiodic point.
    \end{itemize}
    Dynabase aims to collect, for both polynomial and rational maps,
    \begin{enumerate}
        \item The full set of distinct directed graphs for each $d$ and $D$.
        \item The longest possible rational periodic cycle for each $d$ and $D$.
        \item The longest possible rational tail for each $d$ and $D$.
    \end{enumerate}
    An algorithm to compute the full set of rational preperiodic points was introduced in \cite{Hutz12} and is implemented in SageMath.

\subsection{Postcritically Finite Maps} \label{sect_pcf}
    By Thurston's rigidity theorem, apart from the flexible Latt\`es maps, postcritically finite maps do not occur in nontrivial families, and Benedetto, Ingram, Jones, and Levy proved that the postcritically finite conjugacy classes (other than the flexible Latt\`es maps) form a set of bounded height in $\calM_d(\Qbar)$ \cite{Benedetto5}. Consequently, for each $d$ and $D$ there are only finitely many such conjugacy classes with a representative defined over a field $K$ with $[K:\Q] \leq D$. There are infinitely many over $\Qbar$; for example, $z^2+c$ is postcritically finite for infinitely many $c \in \Qbar$. The full set has been computed in a handful of special cases, with some care needed regarding the difference between finding $K$-rational representatives with a parameterization versus finding all conjugacy classes of postcritically finite maps that have some representative defined over $K$. Dynabase aims to collect the full set of postcritically finite conjugacy classes for each $d$ and $D$.

    Postcritically finite maps play the role of special points in the moduli space, analogous to CM points on modular curves \cite{Baker7}. Over $\Q$, the complete lists are known for polynomials of degree $2$ and $3$ \cite{Anderson, Ingram2} and of degree $4$ \cite{Fraser}, and for quadratic rational maps \cite{Lukas}.

    An algorithm to check whether a given map is postcritically finite is implemented in SageMath. It proceeds by determining the set of critical points, then checking if each is preperiodic. Consequently, for high degree maps, this can be quite expensive as the field of definition for the critical points grows. This is an area where some work on the computational algorithms could help speed up the computation of the data for the collected maps.

\subsection{Automorphism Groups} \label{sect_aut}

    As noted above, the automorphism group is finite, since each automorphism must permute the set of periodic points of each period. This also provides an algorithm for computing the full automorphism group of any given map. This method can be further refined by considering multipliers.
    \begin{defn}
        Let $f \in \Hom_d$ and $P \in \Per_n(f)$. Choose coordinates so that $P$ is not the point at infinity, and let $\bar{f}$ be the corresponding dehomogenization of $f$. The \emph{multiplier} of $P$ (as a fixed point of $f^n$) is
        \begin{equation*}
            \lambda_P = (\bar{f}^n)'(P),
        \end{equation*}
        where $'$ is the usual derivative. It does not depend on the choice of coordinates. If the orbit of $P$ avoids the point at infinity, then $\lambda_P = \prod_{i=0}^{n-1} \bar{f}'(\bar{f}^i(P))$.
    \end{defn}
    It can be further shown that automorphisms preserve the multipliers, so that the number of permutations of $\Per_n(f)$ that must be checked can be significantly decreased, leading to a more efficient algorithm. A similar algorithm is given in Faber-Manes-Viray \cite{FMV} and implemented in SageMath.

    Dynabase aims to collect the following information.
    \begin{enumerate}
        \item For each finite group $G \subset \PGL_2$, the map of smallest degree and smallest field of definition for which $G$ is realized as $\Aut(f)$.
        \item For a given degree $d$ and field degree $D$, all the automorphism groups that are realizable.
    \end{enumerate}
    De Faria and Hutz \cite{Hutz16} provide additional tools to construct maps with a given automorphism group.

    The maps with a nontrivial automorphism group form proper subvarieties of $\calM_d$, the automorphism loci \cite{Miasnikov}, which have been determined completely for $d=2$ \cite{Milnor} and for $d = 3, 4$ \cite{Hutz25}. For polynomials, the possible automorphism groups are classified by Fujimura-Nishizawa \cite{Fujimura3}.

\subsection{Small Height Ratio} \label{sect_height}

    Lang's conjecture for elliptic curves states that for any number field $K$ there is a constant $C>0$, depending only on $K$, such that for every elliptic curve $E$ defined over $K$ and every non-torsion point $P \in E(K)$,
        \begin{equation*}
            \hhat(P) \geq C \log(\mathcal{N}_K(\mathcal{D}_E)),
        \end{equation*}
        where $\mathcal{D}_E$ is the minimal discriminant of the elliptic curve $E$ and $\mathcal{N}_K$ is its norm. In the dynamical systems case, some additional care needs to be taken: the canonical height is invariant under conjugation, $\hhat_{f^{\alpha}}(\alpha^{-1}(P)) = \hhat_f(P)$, while the height of the coefficients of $f$ is not, so the lower bound must be in terms of a height on the moduli space.
        A version of Lang's conjecture for dynamical systems is as follows:
     \begin{conj}[\textup{\cite[Conjecture 4.98]{Silverman10}}] \label{conj-lang}
        Let $K$ be a number field and $d \geq 2$. There is a constant $M(d,K)>0$ depending only on $d$ and $K$ so that for any $f \in \Hom_d(K)$ and any wandering point $Q \in \P^1(K)$ we have
        \begin{equation} \label{eq1}
            \hhat_f(Q) \geq M(d,K) h_{\mathcal{M}}(f)
        \end{equation}
        where $h_{\mathcal{M}}(\cdot)$ is a height on the moduli space of dynamical systems of degree $d$.
    \end{conj}
    One way to determine a height on moduli space is to choose an embedding of $\calM_d$ into an affine or projective space. For example, Milnor proved that $\calM_2(\C) \cong \A^2(\C)$ \cite{Milnor}. Hence, one could take $h([f])$ to be the height of the associated point in $\A^2$. To examine this conjecture computationally, a consistent moduli height needs to be chosen for all $d$ and $D$. Dynabase chooses a height associated to the multipliers. More specifically, define the $n$-multiplier spectrum as the multiset (counted with multiplicity)
    \begin{equation*}
        \Lambda_n = \{ \lambda_P \col P \in \Per_n(f)\},
    \end{equation*}
    where $\lambda_P$ is the multiplier of $P$ as a fixed point of $f^n$. Then define
    \begin{equation*}
        \bsigma_n = (\sigma_{n,1}, \ldots, \sigma_{n,d^n+1}),
    \end{equation*}
    where $\sigma_{n,i}$ is the $i$-th elementary symmetric function evaluated on $\Lambda_n$. These $\bsigma_n$ are invariant under conjugation, so they depend only on the conjugacy class. McMullen proved that, apart from the flexible Latt\`es maps, taking these $\bsigma_n$ for large enough $n$ produces a finite-to-one mapping from $\calM_d$ to an affine space \cite{McMullen2}. This has been refined in various ways.
    \begin{itemize}
        \item Fujimura proved that $\bsigma_1$ is enough to determine finitely many conjugacy classes of polynomials ($(d-2)!$ classes) \cite{Fujimura2}. Huguin proved that $(\bsigma_1,\bsigma_2)$ is enough to be generically one-to-one for polynomials; in fact, the morphism on the moduli space of polynomials is finite and birational onto its image \cite{Huguin2}.
        \item Ji and Xie conjectured that $(\bsigma_1,\bsigma_2)$ is enough to be generically one-to-one for all of $\calM_d$ \cite{Ji}, and Zhang recently proved this conjecture: the morphism is birational onto the closure of its image \cite{Zhang5}.
    \end{itemize}
    Dynabase uses the simpler-to-compute height
    \begin{equation*}
        h_{\calM}(f) = \max_i h(\sigma_{1,i}),
    \end{equation*}
    the largest absolute logarithmic height of the first multiplier invariants, and collects the smallest known examples of the ratio
    \begin{equation*}
        \frac{\hhat_f(P)}{h_{\mathcal{M}}(f)}.
    \end{equation*}
    Note that Looper, assuming the $abcd$-conjecture, proved that a uniform lower bound exists on this ratio for polynomial maps, but the bound is not effective \cite{Looper3}. Searches for small ratios were carried out for quadratic polynomials and quadratic rational maps \cite{Benedetto8} and for cubic polynomials \cite{Benedetto4}, and a genetic algorithm extended them to polynomials of degree up to $12$ and rational maps of degree up to $5$ \cite{Hutz29}.

\section{Scope and Inclusion} \label{sect_scope}
Significant effort was made to search the literature for all known examples related to these four problems. A map is stored once for each conjugacy class over its field of definition $K$: two maps that are conjugate by some $\alpha \in \PGL_2(K)$ form a single entry, while rational twists over $K$ are separate entries, since they can differ for questions involving rationality. The maps are subdivided along several categories:
    \begin{itemize}
        \item The dimension $N$ of the domain $\P^N$. Currently $N=1$.
        \item The degree $d$ of the map.
        \item The type of the map. A map is \emph{polynomial} if it has a totally ramified fixed point over $\Kbar$. Every other map is \emph{non-polynomial}.
        \item The field of definition. Tables are organized first by the degree $D = [K:\Q]$ and then by the field $K$ itself, identified by its LMFDB label \cite{LMFDB}.
    \end{itemize}

\subsection{Distinct examples}
    The goal is to record every extreme behavior that is known to occur, not every map that exhibits it. What is considered a new behavior depends on the problem.
    \begin{itemize}
        \item \emph{Rational preperiodic points.} The basic object is a pair $(f,L)$ of a map $f$ defined over $K$ and a number field $L \supseteq K$, together with the directed graph of the $L$-rational preperiodic points of $f$. The goal is to collect all possible distinct graph structures for a given subproblem. Graphs are compared up to isomorphism as unlabelled directed graphs.
            
            Note that a map can contribute data over several fields; for example, many of the quadratic field examples of \cite{DFK} are maps $z^2+c$ with $c \in \Q$ which have different graph structures as the field $L$ changes. A graph that is realized over multiple fields of the same degree is kept as data for each such field. 

        \item \emph{Postcritically finite maps.} The goal is a complete list of conjugacy classes for each subproblem. The summary tables count conjugacy classes over $\Qbar$ and also report the number of maps when twists make the two counts differ.
            
        \item \emph{Automorphism groups.} For each subproblem, one map is presented which realizes each automorphism group.
            
        \item \emph{Small height ratio.} The map and point which realize the smallest known canonical height ratio are presented for each subproblem.
    \end{itemize}
    A list of citations is attached to each map. These citations include any source which has studied that map (including possibly in a different but $K$-conjugate form). In this way, each map is credited to every source in which it appears.

\subsection{Sources}
    Maps are mainly collected from published articles, but preprints and personal communications are also accepted. Each map is entered in a data file together with the citation labels of its sources, so that every displayed example links to where it appears in the literature. An extreme example that can be verified by computation does not need a published source to be included. However, a subproblem is only marked ``proven'' if the complete answer has a published proof. New data can be submitted through a form on the website or by a pull request on the public GitHub project.

\subsection{Status of a subproblem} \label{sect_status}
    For each problem, degree, type and field degree, the website rates the status of the subproblem as one of the following.
    \begin{description}
        \item[proven] The complete answer is a theorem in a published source, and the database contains the complete list of examples.
            
        \item[conjectural] A complete answer has been conjectured, in a published source.
        \item[experimental] The examples come from computer searches, with no conjectured complete answer.
            
        \item[uncomputed] The complete answer is a theorem, but the database does not yet contain its examples.
            
        \item[open] Little or no systematic data is available.
    \end{description}
    The ratings are curated by hand from the literature and cite their sources. The record holders, in contrast, are computed from the data, so they update automatically as new maps are added.

\section{Computation} \label{sect_computation}
    Each map is entered as a pair of homogeneous polynomials over a number field, processed in SageMath \cite{sage2}, and the results are stored in a PostgreSQL database \cite{PostgreSQL}. This section describes the normalizations that make the entries comparable and the identification and inclusion criteria. Note that a few maps were computed via specialized code, notably the maps from \cite{JDoyle4}. These are all contained in the \texttt{add\_functions\_special\_cases\_dim\_1.py} file. There are also some minor improvements such as in automorphism group calculations and postcritically finite verification that have not yet made it into the general SageMath codebase.

\subsection{Normalizing the field}
    The defining polynomial of the field is first reduced with PARI's \texttt{polredabs}, and the field is identified by its LMFDB label \cite{LMFDB}. A field that is not in the LMFDB is labelled by its reduced defining polynomial instead. All coordinates, of maps and of points, are expressed in terms of a root of the reduced polynomial, so a given field always has the same presentation in the database.

\subsection{Normalizing the map}
    Besides the representation as entered (the \emph{original model}), Dynabase computes several normal forms when they exist:
    \begin{itemize}
        \item the monic centered model of a polynomial, $z^d + a_{d-2}z^{d-2} + \cdots + a_0$;
        \item the reduced model of \cite{Hutz20} for maps over $\Q$, which minimizes the size of the coefficients within the $\mathrm{SL}_2(\Z)$-orbit;
        \item whether the map is a Chebyshev map or a Newton map, with the corresponding polynomial.
    \end{itemize}
    One of the models is chosen as the \emph{display model} shown on the website.

\subsection{Identifying maps} \label{sect_identify}
    Similar to the automorphism group computations, Faber-Manes-Viray \cite{FMV} provide algorithms to determine if two maps are conjugate. So, determining if a map is already in the database comes down to checking this conjugacy. However, this can be quickly optimized by checking the following.
    \begin{enumerate}
        \item If one of the stored models (original, monic centered or reduced) of a candidate has exactly the same coefficients over the same field, the map is already present.
        \item Compute $\bsigma_1$ and its hash (Section \ref{sect_labels}). If no stored map of the same degree has the same hash, the map is new.
        \item Compute $\bsigma_2$ and its hash (Section \ref{sect_labels}). If no stored map of the same degree has the same pair of hashes, the map is new.
    \end{enumerate}
    Finally, the set of maps with the same pair $(\bsigma_1,\bsigma_2)$ can be checked via the FMV algorithm.
    
    If the map is found to already be in the database but has a new citation, the citations are updated and the map is kept in the data files.

\section{Labels} \label{sect_labels}
    An important aspect of Cremona's original elliptic curve tables \cite{Cremona} was allowing for the normalization and labelling of elliptic curves. The LMFDB expanded this labelling system and also includes a labelling system for number fields. Similarly, Dynabase proposes a labelling system for algebraic dynamical systems. Based on the theorem of Zhang, conjectured by Ji and Xie, that $(\bsigma_1, \bsigma_2)$ is sufficient for a generically one-to-one identification of conjugacy classes \cite{Ji,Zhang5}, we propose the following label
    \begin{equation*}
        \texttt{dim.sigma1.sigma2.ord}
    \end{equation*}
    where
    \begin{enumerate}
        \item \texttt{dim} is $N$, where the domain is $\P^N$.
        \item \texttt{sigma1} is a hash of the fixed point multiplier invariants $\bsigma_1$.
        \item \texttt{sigma2} is a hash of the $2$-periodic point multiplier invariants $\bsigma_2$.
        \item \texttt{ord} is a positive integer to handle rational twists and other isospectral situations.
    \end{enumerate}
    The hashes are computed as follows:
    \begin{enumerate}
        \item the SageMath string representations of the invariants in $\bsigma_n$ are concatenated (separated by `,') and the resulting string is UTF-8 encoded;
        \item the encoded string is hashed with SHAKE-256 to a $4$-byte digest, written as $8$ hex characters.
    \end{enumerate}
    Note that the normalization of the field ensures that the representations of the $\bsigma_n$ do not cause trivial collisions.
    \begin{exmp}
        A label for $f(z) = \frac{2}{z^2}$ is given as
        \begin{equation*}
            \texttt{1.db03b80b.9e054174.1}
        \end{equation*} 
\begin{python}
sage: import hashlib
sage: def sigma_hash(F, n):
....:     s = ','.join(str(t) for t in F.sigma_invariants(n))
....:     return hashlib.shake_256(s.encode('utf-8')).hexdigest(4)
sage: P.<x,y> = ProjectiveSpace(QQ, 1)
sage: g = DynamicalSystem([2*y^2, x^2])
sage: sigma_hash(g, 1), sigma_hash(g, 2)
('db03b80b', '9e054174')
\end{python}
    \end{exmp}
    \begin{rem}
        The separators are needed for the hash since it is not too hard to create a collision without the separator.
        \begin{align*}
            f &= (2x^2 - 20xy - 20y^2 : -x^2 + 24xy - 222y^2) \qquad \bsigma_1 = (22,202,20)\\
            g &= (2x^2 - 220xy - 220y^2 : -x^2 + 224xy - 220y^2) \qquad \bsigma_1 =  (222, 0, 220)
        \end{align*}    
    \end{rem}
    
    The ordinal \texttt{ord} distinguishes maps with the same degree and the same pair of hashes: rational twists, and also any non-conjugate maps with the same multipliers, such as the flexible Latt\`es maps. It is assigned in the order in which the maps were added and is then recorded in a label file kept with the data files, so that rebuilding the database, even with the data files loaded in a different order, gives each map the same label. As with LMFDB labels, a label is permanent once assigned.

\section{The Data} \label{sect_snapshot}
    This section describes the contents of version 1.0.0 of the data. This data was computed on a laptop running an Intel i7-1360P with 32 GB of RAM with 8 parallel processes with SageMath running under the Windows Subsystem for Linux. Note that only about half of the RAM was available to WSL. The software versions are listed in Table \ref{table_software}.
    \begin{table}[ht]
    \small
    \begin{tabular}{|l|l|}
        \hline
        Software & Version\\
        \hline
        WSL - Ubuntu & 22.04\\
        \hline
        SageMath & 10.4\\
        \hline
        PostgreSQL & 14.24\\
        \hline
        psycopg2 & 2.9.9\\
        \hline
        lmfdb-lite & 0.1\\
        \hline
    \end{tabular}
    \caption{Software versions.} \label{table_software}
    \end{table}

    It took approximately 2 days of wall time to compute this dataset with a 7200 second timeout. All properties related to the 4 main problems are complete, but a number of additional properties did not complete. The main missing data values are alternate models (reduced or monic-centered) for a map and the field of definition of the critical points. The vast majority of the properties successfully compute with a 300s timeout for those wanting to recompute a nearly full dataset in a shorter time frame. An export of the dataset is available both on GitHub\footnote{\url{https://github.com/bhutz/DynaBase}} and Zenodo\footnote{DOI: \url{10.5281/zenodo.24377332}}. This article can be cited for references to the general Dynabase project. The GitHub project is licensed under GPL version 3 or later and the data under CC BY 4.0.

\subsection{Summary}
    The database contains 592 maps from 28 sources, 556 of them over $\Q$ and 36 over 17 other fields: 29 maps over 13 quadratic fields and 7 maps over fields of degree 4, 6, 8 and 12. Table \ref{table_counts} breaks the maps down by degree and type. There are 633 pairs of a map and a field with their rational preperiodic points, 69 of them over quadratic fields, which realize 366 distinct graph structures. Of the maps, 73 are postcritically finite, 44 have a smallest height ratio, and 192 have a nontrivial automorphism group, realizing 19 distinct groups.

    \begin{table}[ht]
    \small
    \begin{tabular}{|l|*{14}{c|}}
        \hline
        $d$ & 2 & 3 & 4 & 5 & 6 & 7 & 8 & 9 & 10 & 11 & 12 & 13 & 14 & 15\\
        \hline
        Polynomial & 44 & 123 & 33 & 18 & 30 & 24 & 22 & 18 & 17 & 13 & 12 & 8 & 1 & 1\\
        \hline  
        Non-polynomial & 72 & 72 & 76 & 5 & 1 & 1 & 0 & 0 & 0 & 1 & 0 & 0 & 0 & 0\\
        \hline
    \end{tabular}
    \caption{The number of maps by degree and type.} \label{table_counts}
    \end{table}

    Table \ref{table_sources} lists the sources, in the order of the bibliography, with the number of maps that cite each. A map can cite several sources.

    \begin{table}[ht]
    \small
    \begin{tabular}{|l|r||l|r||l|r|}
        \hline
        Source & Maps & Source & Maps & Source & Maps\\
        \hline
        \cite{Anderson} & 15 & \cite{JDoyle4} & 12 & \cite{Miasnikov} & 1\\
        \hline
        \cite{Benedetto8} & 20 & \cite{Doyle-Hyde} & 15 & \cite{Milnor} & 1\\
        \hline
        \cite{Benedetto4} & 105 & \cite{Fraser} & 16 & \cite{Morton2} & 2\\
        \hline
        \cite{Bilgili} & 4 & \cite{Fujimura3} & 23 & \cite{Poonen} & 12\\
        \hline
        \cite{Blanc} & 1 & \cite{Hutz25} & 118 & \cite{Stoll3} & 1\\
        \hline
        \cite{Canci6} & 24 & \cite{Hutz29} & 204 & \cite{Vishkautsan} & 6\\
        \hline
        \cite{Hutz16} & 15 & \cite{Hutz5} & 2 & \cite{Vishkautsan2} & 1\\
        \hline
        \cite{JDoyle} & 38 & \cite{Ingram2} & 9 & \cite{Walde} & 4\\
        \hline
        \cite{JDoyle3} & 13 & \cite{Lukas} & 12 & & \\
        \hline
        \cite{DFK} & 38 & \cite{Manes2} & 10 & & \\
        \hline
    \end{tabular}
    \caption{Sources of the data and the number of maps that cite each.} \label{table_sources}
    \end{table}

\subsection{Rational preperiodic points}
    Table \ref{table_preper} gives the records over $\Q$: for each degree and type, the most rational preperiodic points, the longest rational cycle and the longest rational tail, together with the number of distinct graph structures in the database. Points are counted in $\P^1(\Q)$, including the point at infinity. The non-polynomial maps of degree at least 5 in the database are only the examples of automorphism groups, so their records say little. Over quadratic fields, the database has quadratic polynomials realizing 46 graph structures, with a record of 15 points, a cycle of length 6 and a tail of length 4 \cite{Walde,Poonen,Stoll3,Hutz5,DFK,JDoyle,JDoyle3}, and two polynomials of degrees 4 and 6 over $\Q(\sqrt{-3})$ \cite{JDoyle4}.
    
    Note that the code attempts to compute all properties for every map. So, until they are superseded by better data, some of the subproblems have a small amount of superficial data appearing.

    \begin{table}[ht]
    \small
    \setlength{\tabcolsep}{4pt}
    \begin{tabular}{|c|c|c|c|c|c||c|c|c|c|c|}
        \hline
        & \multicolumn{5}{c||}{Polynomial} & \multicolumn{5}{c|}{Non-polynomial}\\
        \hline
        $d$ & Graphs & Points & Cycle & Tail & Sources & Graphs & Points & Cycle & Tail & Sources\\
        \hline
        2 & 12 & 9 & 3 & 2 & \cite{JDoyle,JDoyle3,DFK,Ingram2}
        & 57 & 16 & 7 & 6 & \cite{Benedetto8,Vishkautsan2}\\
        &&&&&\cite{Lukas,Poonen,Walde}&&&&&\\
        \hline
        3 & 104 & 12 & 5 & 4 & \cite{Benedetto4,Doyle-Hyde,Hutz29} & 59 & 13 & 8 & 8 & \cite{Hutz29}\\
        \hline
        4 & 17 & 11 & 6 & 5 & \cite{Doyle-Hyde,Hutz29} & 68 & 12 & 10 & 9 & \cite{Hutz25,Hutz29}\\
        \hline
        5 & 14 & 14 & 8 & 6 & \cite{Doyle-Hyde,Hutz29} & 3 & 4 & 2 & 0 & \cite{Hutz16}\\
        \hline
        6 & 17 & 15 & 8 & 7 & \cite{Doyle-Hyde,Hutz29} & 1 & 4 & 2 & 1 & \cite{Hutz16}\\
        \hline
        7 & 17 & 16 & 10 & 9 & \cite{Doyle-Hyde,Hutz29} & 1 & 4 & 2 & 0 & \cite{Hutz16}\\
        \hline
        8 & 11 & 17 & 11 & 9 & \cite{Doyle-Hyde,Hutz29} & & & & &\\
        \hline
        9 & 7 & 20 & 11 & 10 & \cite{Doyle-Hyde,Hutz29} & & & & &\\
        \hline
        10 & 7 & 21 & 12 & 11 & \cite{Doyle-Hyde,Hutz29} & & & & &\\
        \hline
        11 & 4 & 20 & 13 & 12 & \cite{Doyle-Hyde,Hutz29} & 1 & 2 & 1 & 0 & \cite{Hutz16}\\
        \hline
        12 & 4 & 21 & 14 & 13 & \cite{Doyle-Hyde,Hutz29} & & & & &\\
        \hline
        13 & 3 & 22 & 15 & 14 & \cite{Doyle-Hyde,Hutz29} & & & & &\\
        \hline
        14 & 1 & 23 & 1 & 4 & \cite{Doyle-Hyde} & & & & &\\
        \hline
        15 & 1 & 24 & 4 & 7 & \cite{Doyle-Hyde} & & & & &\\
        \hline
    \end{tabular}
    \caption{Records for rational preperiodic points over $\Q$: the number of distinct graph structures, the most rational preperiodic points, the longest cycle and the longest tail, with the sources of the record holders.} \label{table_preper}
    \end{table}

\subsection{Postcritically finite maps}
    Table \ref{table_pcf} gives the number of postcritically finite conjugacy classes over $\Q$. The classifications are complete for polynomials of degree $2$ (3 classes \cite{Ingram2}), $3$ (15 classes \cite{Anderson,Ingram2}) and $4$ (16 classes \cite{Fraser}), and for quadratic rational maps (12 classes in all, of which 9 are non-polynomial \cite{Lukas}). The cubic polynomials over $\Q$ include three pairs of twists, which is why the 15 classes correspond to 18 maps. In degree $5$ and higher, the database has only isolated examples. There are also postcritically finite polynomials of degrees 2, 6, 7 and 11 over fields of degree 2, 4, 6 and 12.

    \begin{table}[ht]
    \small
    \begin{tabular}{|l|*{10}{c|}}
        \hline
        $d$ & 2 & 3 & 4 & 5 & 6 & 7 & 8 & 9 & 10 & 11\\
        \hline
        Polynomial & 3 & 15 & 16 & 1 & 1 & 1 & 1 & 1 & 1 & 0\\
        \hline
        Non-polynomial & 9 & 2 & 1 & 2 & 1 & 1 & 0 & 0 & 0 & 1\\
        \hline
    \end{tabular}
    \caption{Postcritically finite conjugacy classes over $\Qbar$ with a representative over $\Q$ in the database. Only degrees $2$ to $4$ for polynomials and degree $2$ for non-polynomial maps are complete.} \label{table_pcf}
    \end{table}

\subsection{Automorphism groups}
    Table \ref{table_aut} lists, for each group realized in the database, the smallest degree of a map over $\Q$ with exactly that automorphism group. We write $D_n$ for the dihedral group of order $2n$. For polynomials, the column is the classification of \cite{Fujimura3}: a polynomial conjugate to $z^d$ has automorphism group $D_{d-1}$, and any other polynomial in monic centered form (no $z^{d-1}$ term) has the cyclic group $C_m$ with $m = \gcd(e-1)$ over the exponents $e$ of its terms. The database contains an example of every group possible for polynomials in each degree up to $10$. For non-polynomial maps over $\Q$, the possible groups are known in degree $2$ \cite{Milnor} and in degrees $3$ and $4$ \cite{Hutz25}, and the database contains an example of each.

    \begin{table}[ht]
    \small
    \begin{tabular}{|c|c|c|c||c|c|c|c|}
        \hline
        Group & Order & Poly & Non-poly & Group & Order & Poly & Non-poly\\
        \hline
        $C_2$ & 2 & 2 & 2 & $C_9$ & 9 & 10 & \\
        \hline
        $C_3$ & 3 & 4 & 3 & $D_5$ & 10 & 6 & 4\\
        \hline
        $C_2 \times C_2$ & 4 & 3 & 3 & $A_4$ & 12 & -- & 3\\
        \hline
        $C_4$ & 4 & 5 & 4 & $D_6$ & 12 & 7 & 5\\
        \hline
        $C_5$ & 5 & 6 &  & $D_7$ & 14 & 8 & 6\\
        \hline
        $C_6$ & 6 & 7 &  & $D_8$ & 16 & 9 & 7\\
        \hline
        $S_3$ & 6 & 4 & 2 & $D_9$ & 18 & 10 & \\
        \hline
        $C_7$ & 7 & 8 &  & $S_4$ & 24 & -- & 5\\
        \hline
        $C_8$ & 8 & 9 &  & $A_5$ & 60 & -- & 11\\
        \hline
        $D_4$ & 8 & 5 & 3 & & & &\\
        \hline
    \end{tabular}
    \caption{The smallest degree of a map over $\Q$ in the database with exactly the given automorphism group. A dash means that no polynomial has this group; an empty cell means that the database has no example.} \label{table_aut}
    \end{table}

\subsection{Small height ratios}
    Table \ref{table_height} gives the smallest height ratio over $\Q$ in each degree and type. All of the records come from \cite{Hutz29}. The degree 2 and 3 polynomial records and the degree 2 non-polynomial record equal the smallest values found earlier in \cite{Benedetto8,Benedetto4}.

    \begin{table}[ht]
    \small
    \begin{tabular}{|c|c|c||c|c|c|}
        \hline
        $d$ & Polynomial & Non-polynomial & $d$ & Polynomial & Non-polynomial\\
        \hline
        2 & $6.6041 \cdot 10^{-3}$ & $4.6574 \cdot 10^{-4}$ & 7 & $1.0564 \cdot 10^{-10}$ & \\
        \hline
        3 & $9.2099 \cdot 10^{-5}$ & $3.0790 \cdot 10^{-6}$ & 8 & $8.7710 \cdot 10^{-13}$ & \\
        \hline
        4 & $2.9015 \cdot 10^{-6}$ & $2.1843 \cdot 10^{-8}$ & 9 & $1.7173 \cdot 10^{-14}$ & \\
        \hline
        5 & $9.1520 \cdot 10^{-9}$ & $3.6941 \cdot 10^{-10}$ & 10 & $9.9368 \cdot 10^{-16}$ & \\
        \hline
        6 & $1.8372 \cdot 10^{-9}$ & & 11 & $7.3446 \cdot 10^{-17}$ & \\
        \hline
         & & & 12 & $5.9715 \cdot 10^{-19}$ & \\
        \hline
    \end{tabular}
    \caption{The smallest height ratios $\hhat_f(P)/h_{\calM}(f)$ over $\Q$ in the database.} \label{table_height}
    \end{table}

\subsection{Status}
    Table \ref{table_status} summarizes the status of each subproblem over $\Q$, as defined in Section \ref{sect_status}.

    \begin{table}[ht]
    \small
    \begin{tabular}{|c|c|c|c|c|c|c|c|c|}
        \hline
        & \multicolumn{2}{c|}{Preperiodic} & \multicolumn{2}{c|}{PCF} & \multicolumn{2}{c|}{Automorphisms} & \multicolumn{2}{c|}{Small height}\\
        \hline
        $d$ & Poly & Non-poly & Poly & Non-poly & Poly & Non-poly & Poly & Non-poly\\
        \hline
        2 & C & E & P & P & P & P & E & E\\
        \hline
        3 & E & O & P & O & P & P & E & E\\
        \hline
        4 & E & O & P & O & P & P & E & E\\
        \hline
        5 & E & O & O & O & P & O & E & E\\
        \hline
        6--10 & E & O & O & O & P & O & E & O\\
        \hline
        11, 12 & E & O & O & O & U & O & E & O\\
        \hline
        13 & E & O & O & O & U & O & O & O\\
        \hline
        14, 15 & O & O & O & O & U & O & O & O\\
        \hline
    \end{tabular}
    \caption{Status of the subproblems over $\Q$: proven (P), conjectural (C), experimental (E), uncomputed (U) and open (O). The three rational preperiodic subproblems (most points, longest cycle, longest tail) have the same status.} \label{table_status}
    \end{table}

\subsection{Data versions} \label{sect_versions}
    Versions of the data are numbered \texttt{major.minor.patch}. The major number increases with a change to the database schema, the minor number when data is added or newly computed, and the patch number with a correction to existing values. Labels are stable from version 1.0.0 on. The label file records the label of every map; a map being loaded is assigned its recorded label. This is to prevent a reordering of the input data resulting in a different ordinal assigned to twists.

    Each version is a tagged commit of the GitHub repository, together with the version of SageMath used. This data is also packaged and added to Zenodo. Since a computation that reaches a time limit leaves a property unset, a rebuild is not guaranteed to reproduce a version exactly, so the tagged commit also contains each table as sorted JSON lines with its internal identifiers replaced by labels, as well as the label file. The website shows the current version at the bottom of each page and lists all versions in a changelog, and each SageMath export records the version it came from.

\section{Corrections to the Literature} \label{sect_corrections}
    Compiling the database required recomputing every example, which uncovered the following errors in the sources.
    \begin{itemize}
        \item It was falsely claimed in \cite{Hutz16} that no map defined over $\Q$ has automorphism group $A_4$ and later corrected in \cite{Hutz25} with the map $(z^3-3)/(-3z^2)$.
        \item In \cite[Figure 1]{Manes2}, example (j) is labelled $(k,b) = (4/3, 2/15)$. The value should be $b = -2/15$, and it is correct in the text of that paper, but not in the figure.
        \item In \cite[Theorem 6.9]{Hutz25}, the structure $G_1$ of the $C_3$ family is given by $a = 1/(1-t^3)$ with $t = 1/(xy-xz)$ for the points $(x:y:z) = n(1:0:1)$, $n \geq 2$, on $y^2 - y = x^3 - 1$. It should be $t = (yz-z^2)/(xy-xz)$, as follows from the map $\psi$ in the proof of the theorem. For example, $n=2$ gives the point $(7:19:1)$, $t = 1/7$ and $a = 343/342$.
    \end{itemize}

\section{Using Dynabase} \label{sect_using}
    The data is published as a set of static web pages at \url{https://dynabase.org}, generated from the database. The pages include
    \begin{itemize}
        \item an introduction and background on the problems;
        \item a summary of the included data, with the sources of every map;
        \item data pages for every problem, degree, type and field, with one table each;
        \item a summary of the extreme examples, with the record holder of every subproblem;
        \item the status charts of every problem, with their citations;
        \item a form for comments and submissions.
    \end{itemize}
    A filter bar selects the dimension, degree, type, field degree and field of a data page. Each map is displayed with its label, its sources, and a link to an image of its graph of rational preperiodic points or its critical portrait. Every table can be exported as a file that, when loaded in SageMath, defines a list \texttt{dynabase\_systems} of pairs consisting of a map and its label. In the export of a small height table, each map also comes with its point of smallest height ratio.

    For example, the following session loads the degree $2$ polynomials over $\Q$ from the rational preperiodic points page and recomputes the preperiodic points of $z^2 - 29/16$.
\begin{python}
sage: # from https://dynabase.org/data/1/2/poly/1/preper.sage
sage: load("preper.sage")
sage: f, label = dynabase_systems[1]
sage: f, label
(Dynamical System of Projective Space of dimension 1 over Rational Field
   Defn: Defined on coordinates by sending (x : y) to
         (x^2 - 29/16*y^2 : y^2), '1.2799614f.68d04ee2.1')
sage: f.rational_preperiodic_graph().num_verts()
9
\end{python}
    % TODO: add a screenshot of a data page and one graph image as figures.

\section{Architecture, Reproducibility, and Contributing} \label{sect_architecture}
    Dynabase is developed on GitHub under the GNU General Public License, version 3 or later. The repository has three parts.
    \begin{itemize}
        \item \emph{The data files.} Each source has a section in a SageMath script that constructs its maps and adds each one with its citation labels; the bibliography is a separate script. A driver script loads the bibliography, fields, families and maps in order. A map already in the database only has its citations updated.
        \item \emph{The computation code.} SageMath functions normalize fields and maps, identify maps, and compute the properties of Section \ref{sect_computation}. Field labels are looked up in the LMFDB through \texttt{lmfdb-lite}.
        \item \emph{The site generator.} A Python program reads the database and writes the static website, including the SageMath exports and the graph images. It requires a local version of SageMath to run, and its output is published by GitHub Pages.
    \end{itemize}
    The database has tables for the maps, for the graph structures (so each graph is stored once and referenced by its rational preperiodic points and critical portraits), for the rational preperiodic points of each pair of a map and a field, for the bibliography, and for one-parameter families of maps. Rebuilding the database requires SageMath, PostgreSQL and access to the LMFDB, after which the database can be populated by running the driver script.

    To contribute, add the maps of a new source as a section of the appropriate data file, with a citation entry for the source, and open a pull request. Alternatively, submit the maps through the form on the website. A contribution is checked by recomputing the example. To cite Dynabase, cite this article and the version of the data used, as listed on the website.

\section{AI Methodology} \label{sect_AI}
    Claude models Opus 5 and 5.5 were used in various stages of this project. It was used mainly as a coding assist tool and was not used to generate new mathematical concepts or ideas. The main uses were to
    \begin{itemize}
        \item Produce the website pages hosted via GitHub Pages and the data packaging for Zenodo.
        \item Improve the parallel code for data generation.
        \item Summarize the data within the computed dataset.
        \item Extract the examples from the published paper sources.
    \end{itemize}
    The GitHub commits that were assisted by Claude are labeled as such in the commit messages. Claude was used in the manuscript only for copy editing and compiling of the data tables for the Data section (Section \ref{sect_snapshot}).

\section{Future Work and Directions}\label{sect_future}
\subsection{Expanding data for the current problems}
Most of the data (proven or experimental) is clustered around low degrees over $\Q$ (or quadratic fields). It would be a good student research project to take almost any of the ``open'' subproblems in rational preperiodic points and experimentally determine an interesting set of examples. The small height examples are more filled in due to \cite{Hutz29}. The PCF examples are more in line with graduate student projects, as their utility comes from having \underline{all} PCF examples for a given subproblem. Finally, the automorphism ``open'' subproblems are the most accessible due to the tools developed in \cite{Hutz16}.

\subsection{Expanding data to other problems}
There are likely other problems in arithmetic or algebraic dynamics that would fit the paradigm of Dynabase. These problems could be added and would need to include building the appropriate tools, finding known examples, and creating the website pages.

\subsection{Higher dimensions}
Expanding the data to include higher dimensions is also a good direction to pursue. Many of the computational tools in SageMath work well in higher dimensions, but fewer examples are known. Similarly, the labelling system would need to be examined, as well as taking into account rational map versus morphism behavior. These details range from small technical details to significant open theoretical work.

\subsection{Finite fields}
    The current problems, data, and tools focus on number fields, but the computational tools in SageMath can also handle the finite field case. Expanding in this direction would need to consider which problems most need to be tracked, build the necessary tools, and examine the literature for known examples related to those problems.

\providecommand\biburl[1]{\texttt{#1}}


\begin{thebibliography}{10}

\bibitem{Anderson}
Jacqueline Anderson, Michelle Manes, and Bella Tobin.
\newblock Cubic post-critically finite polynomials defined over $\mathbb{Q}$.
\newblock In {\em ANTS XIV---Proceedings of the Fourteenth Algorithmic Number
  Theory Symposium}, volume~4 of {\em Open Book Ser.}, pages 23--38. Math. Sci.
  Publ., Berkeley, CA, 2020.

\bibitem{Baker7}
Matthew Baker and Laura DeMarco.
\newblock Special curves and postcritically-finite polynomials.
\newblock {\em Forum of Mathematics, Pi}, 1:1--35, 2013.

\bibitem{Benedetto8}
Robert Benedetto, Ruqian Chen, Trevor Hyde, Yordanka Kovacheva, and Colin
  White.
\newblock Small dynamical heights for quadratic polynomials and rational
  functions.
\newblock {\em Experimental Mathematics}, 23(4):433--447, 2014.

\bibitem{Benedetto4}
Robert Benedetto, Ben Dickman, Sasha Joseph, Ben Krause, Dan Rubin, and Xinwen
  Zhou.
\newblock Computing points of small height for cubic polynomials.
\newblock {\em Involve}, 2:37--64, 2009.

\bibitem{Benedetto5}
Robert Benedetto, Patrick Ingram, Rafe Jones, and Alon Levy.
\newblock Attracting cycles in {$p$}-adic dynamics and height bounds for
  postcritically finite maps.
\newblock {\em Duke Math. J.}, 163(13):2325--2356, 2014.

\bibitem{Bilgili}
Hasan Bilgili and Mohammad Sadek.
\newblock Rational preperiodic points of quadratic rational maps over
  {$\mathbb{Q}$} with nonabelian automorphism groups.
\newblock arXiv:2603.06203, 2026.

\bibitem{Blanc}
J{\'e}r{\'e}my Blanc, Jung~Kyu Canci, and Noam~D. Elkies.
\newblock Moduli spaces of quadratic rational maps with a marked periodic point
  of small order.
\newblock {\em Int. Math. Res. Not. IMRN}, 2015(23):12459--12489, 2015.

\bibitem{Booker}
Andrew~R. Booker, Jeroen Sijsling, Andrew~V. Sutherland, John Voight, and Dan
  Yasaki.
\newblock A database of genus-2 curves over the rational numbers.
\newblock {\em LMS J. Comput. Math.}, 19(suppl. A):235--254, 2016.

\bibitem{Canci6}
Jung~Kyu Canci and Solomon Vishkautsan.
\newblock Quadratic maps with a periodic critical point of period 2.
\newblock {\em International Journal of Number Theory}, 13(6):1393--1417, 2017.

\bibitem{Cremona}
J.~E. Cremona.
\newblock {\em Algorithms for modular elliptic curves}.
\newblock Cambridge University Press, 1992.

\bibitem{Hutz16}
Joao de~Faria and Benjamin Hutz.
\newblock Automorphism groups and invariant theory on {$\mathbb{P}^N$}.
\newblock {\em Journal of Algebra and Its Applications}, 17(9):1850162, 38,
  2018.

\bibitem{JDoyle}
John~R. Doyle.
\newblock Preperiodic points for quadratic polynomials with small cycles over
  quadratic fields.
\newblock {\em Math. Z.}, 289(1-2):729--786, 2018.

\bibitem{JDoyle3}
John~R. Doyle.
\newblock Preperiodic points for quadratic polynomials over cyclotomic
  quadratic fields.
\newblock {\em Acta Arith.}, 196(3):219--268, 2020.

\bibitem{DFK}
John~R. Doyle, Xander Faber, and David Krumm.
\newblock Preperiodic points for quadratic polynomials over quadratic fields.
\newblock {\em New York Journal of Mathematics}, 20:507--605, 2014.

\bibitem{JDoyle4}
John~R. Doyle and Wade Hindes.
\newblock Preperiodic integers for {$x^d+c$} in large degree.
\newblock arXiv:2510.14067, 2025.

\bibitem{Doyle-Hyde}
John~R. Doyle and Trevor Hyde.
\newblock Polynomials with many rational preperiodic points.
\newblock {\em Transactions of the American Mathematical Society},
  378(8):5981--6011, 2025.

\bibitem{Dujella}
Andrej Dujella.
\newblock High rank elliptic curves with prescribed torsion.
\newblock \url{https://web.math.pmf.unizg.hr/~duje/tors/tors.html}, 2026.

\bibitem{FMV}
Xander Faber, Michelle Manes, and Bianca Viray.
\newblock Computing conjugating sets and automorphism groups of rational
  functions.
\newblock {\em Journal of Algebra}, 423:1161--1190, 2015.

\bibitem{Fakhruddin}
Najmuddin Fakhruddin.
\newblock Questions on self-maps of algebraic varieties.
\newblock {\em J. Ramanujan Math. Soc.}, 18(2):109--122, 2003.

\bibitem{FPS}
E.V. Flynn, Bjorn Poonen, and Edward Schaefer.
\newblock Cycles of quadratic polynomials and rational points on a genus 2
  curve.
\newblock {\em Duke Math. J.}, 90:435--463, 1997.

\bibitem{Fraser}
Benjamin Fraser.
\newblock {\em On {PCF} polynomials}.
\newblock PhD thesis, York University, Toronto, 2024.

\bibitem{Fujimura2}
Masayo Fujimura.
\newblock The moduli space of rational maps and surjectivity of multiplier
  representation.
\newblock {\em Comput. Methods Funct. Theory}, 7(2):345--360, 2007.

\bibitem{Fujimura3}
Masayo Fujimura and Kiyoko Nishizawa.
\newblock Moduli spaces and symmetry loci of polynomial maps.
\newblock In {\em Proceedings of the 1997 International Symposium on Symbolic
  and Algebraic Computation}, ISSAC '97, pages 342--348, New York, NY, USA,
  1997. ACM.

\bibitem{Hutz25}
Brandon Gontmacher, Benjamin Hutz, Grayson Jorgenson, Srinjoy Srimani, and
  Simon Xu.
\newblock Automorphism loci for degree 3 and degree 4 endomorphisms of the
  projective line.
\newblock {\em New York Journal of Mathematics}, 27:1613--1702, 2021.

\bibitem{Huguin2}
Valentin Huguin.
\newblock Moduli spaces of polynomial maps and multipliers at small cycles.
\newblock arXiv:2412.19335, 2024.

\bibitem{Hutz12}
Benjamin Hutz.
\newblock Determination of all rational preperiodic points for morphisms of
  {$\mathbb{P}^N$}.
\newblock {\em Mathematics of Computation}, 84(291):289--308, 2015.

\bibitem{Hutz29}
Benjamin Hutz.
\newblock A genetic algorithm for generating extreme examples in arithmetic
  dynamics.
\newblock {\em arXiv:2601.11482}, 2026.

\bibitem{Hutz5}
Benjamin Hutz and Patrick Ingram.
\newblock On {P}oonen's conjecture concerning rational preperiodic points of
  quadratic maps.
\newblock {\em Rocky Mountain Journal of Mathematics}, 43(1):193--204, 2013.

\bibitem{Hutz20}
Benjamin Hutz and Michael Stoll.
\newblock Smallest representatives of {${\rm SL}(2,\mathbb{Z})$}-orbits of
  binary forms and endomorphisms of {$\mathbb{P}^1$}.
\newblock {\em Acta Arithmetica}, 189(3):283--308, 2019.

\bibitem{Ingram2}
Patrick Ingram.
\newblock A finiteness result for post-critically finite polynomials.
\newblock {\em Int. Math. Res. Not.}, 3:524--543, 2012.

\bibitem{Ji}
Zhuchao Ji and Junyi Xie.
\newblock The multiplier spectrum morphism is generically injective.
\newblock {\em arxiv:2309.15382}, 2023.

\bibitem{Jones4}
John~W. Jones and David~P. Roberts.
\newblock A database of number fields.
\newblock {\em LMS J. Comput. Math.}, 17(1):595--618, 2014.

\bibitem{Levy}
Alon Levy.
\newblock The space of morphisms on projective space.
\newblock {\em Acta Arithmetica}, 146:12--31, 2011.

\bibitem{Looper3}
Nicole Looper.
\newblock The uniform boundedness and dynamical {L}ang conjectures for
  polynomials.
\newblock arXiv:2105.05240, 2021.
\newblock To appear in Algebra and Number Theory.

\bibitem{Lukas}
David Lukas, Michelle Manes, and Diane Yap.
\newblock A census of quadratic post-critically finite rational functions
  defined over $\mathbb{Q}$.
\newblock {\em LMS Journal of Computation and Mathematics}, 17(A):314--329,
  2014.

\bibitem{Manes2}
Michelle Manes.
\newblock $\mathbb{Q}$-rational cycles for degree-2 rational maps having an
  automorphism.
\newblock {\em Proc. London Math. Soc.}, 96:669--696, 2008.

\bibitem{McMullen2}
Curtis McMullen.
\newblock Families of rational maps and iterative root-finding algorithms.
\newblock {\em Ann. of Math.}, 125(3):467--493, 1987.

\bibitem{Miasnikov}
Nikita Miasnikov, Brian Stout, and Phillip Williams.
\newblock Automorphism loci for the moduli space of rational maps.
\newblock {\em Acta Arithmetica}, 180(3):267--296, 2017.

\bibitem{Milnor}
J.~Milnor.
\newblock Geometry and dynamics of quadratic rational maps.
\newblock {\em Experiment. Math.}, 2(1):37--83, 1993.

\bibitem{Morton2}
Patrick Morton.
\newblock Arithmetic properties of periodic points of quadratic maps.
\newblock {\em Acta Arith.}, 62:343--372, 1992.

\bibitem{Morton5}
Patrick Morton.
\newblock Arithmetic properties of periodic points of quadratic maps, {II}.
\newblock {\em Acta Arith.}, 87:89--102, 1998.

\bibitem{Silverman7}
Patrick Morton and Joseph~H. Silverman.
\newblock Rational periodic points of rational functions.
\newblock {\em Int. Math. Res. Not.}, 2:97--110, 1994.

\bibitem{Northcott}
D.G. Northcott.
\newblock Periodic points of an algebraic variety.
\newblock {\em Ann. of Math.}, 51:167--177, 1950.

\bibitem{petsche}
Clayton Petsche, Lucien Szpiro, and Michael Tepper.
\newblock Isotriviality is equivalent to potential good reduction for
  endomorphisms of {$\mathbb{P}^N$} over function fields.
\newblock {\em Journal of Algebra}, 322(9):3345--3365, 2009.

\bibitem{Poonen}
Bjorn Poonen.
\newblock The complete classification of rational preperiodic points of
  quadratic polynomials over $\mathbb{Q}$: a refined conjecture.
\newblock {\em Math. Z.}, 228(1):11--29, 1998.

\bibitem{Silverman12}
Joseph~H. Silverman.
\newblock The field of definition for dynamical systems on $\mathbb{P}^{1}$.
\newblock {\em Comp. Math.}, 98:269--304, 1995.

\bibitem{Silverman9}
Joseph~H. Silverman.
\newblock The space of rational maps on {$\mathbb{P}^1$}.
\newblock {\em Duke Math. J.}, 94:41--118, 1998.

\bibitem{Silverman10}
Joseph~H. Silverman.
\newblock {\em The Arithmetic of Dynamical Systems}, volume 241 of {\em
  Graduate Texts in Mathematics}.
\newblock Springer-Verlag, 2007.

\bibitem{Stoll3}
Michael Stoll.
\newblock Rational 6-cycles under iteration of quadratic polynomials.
\newblock {\em London Math. Soc. J. Comput. Math.}, 11:367--380, 2008.

\bibitem{LMFDB}
{The LMFDB Collaboration}.
\newblock The {L}-functions and modular forms database.
\newblock \url{https://www.lmfdb.org}, 2026.

\bibitem{PostgreSQL}
{The PostgreSQL Global Development Group}.
\newblock {\em {PostgreSQL}}, 2026.
\newblock \url{https://www.postgresql.org}.

\bibitem{sage2}
{The Sage Developers}.
\newblock {\em {SageMath}, the {S}age {M}athematics {S}oftware {S}ystem
  ({V}ersion 10.4)}, 2024.
\newblock \url{https://www.sagemath.org}.

\bibitem{Vishkautsan}
Solomon Vishkautsan.
\newblock Quadratic rational functions with a rational periodic critical point
  of period 3.
\newblock {\em J. Th\'{e}or. Nombres Bordeaux}, 31(1):49--79, 2019.
\newblock With an appendix by Michael Stoll.

\bibitem{Vishkautsan2}
Solomon Vishkautsan.
\newblock Personal communication, 2026.

\bibitem{Walde}
Ralph Walde and Paula Russo.
\newblock Rational periodic points of the quadratic function {$Q_c(x) = x^2 +
  c$}.
\newblock {\em Amer. Math. Monthly}, 101(4):318--331, 1994.

\bibitem{Zhang5}
Geng-Rui Zhang.
\newblock Generic reconstruction of rational maps from multipliers of periods
  one and two.
\newblock arXiv:2609.07383, 2026.

\end{thebibliography}
\end{document}